\documentclass[a4paper,11pt]{article}
\usepackage{authblk}
\usepackage[latin1]{inputenc}
\usepackage[pdftex]{graphicx}
\usepackage{subfig}
\usepackage{float}
\usepackage{amsfonts,amssymb,amsmath,amsthm}
\usepackage{pdfsync}
\usepackage{color}
\usepackage{url}
\usepackage{nicefrac}
\usepackage[pdftex]{hyperref}
\usepackage{hyperref}

 \hypersetup{colorlinks=true,citecolor=blue, filecolor=red, linkcolor=blue, urlcolor=blue}

\newtheorem{theorem}{Theorem}[section] 
\newtheorem{proposition}[theorem]{Proposition}
\newtheorem{lemma} [theorem] {Lemma}

\newtheorem{remark}[theorem]{Remark}
\newtheorem{corollary}[theorem]{Corollary}

\newcommand \Tr{{\rm Tr\,}}

\makeatletter
\@addtoreset{equation}{section}

\makeatother

\numberwithin{equation}{section}

\newcommand{\R}{\mathbb{R}}
\newcommand{\N}{\mathbb{N}}

\newcommand{\hexagon}{{hex}}

\newcommand{\be}{\begin{equation}}
\newcommand{\ee}{\end{equation}}
\newcommand{\ben}{\begin{equation*}}
\newcommand{\een}{\end{equation*}}
\newcommand{\ba}{\begin{eqnarray}}
\newcommand{\ea}{\end{eqnarray}}
\newcommand{\ban}{\begin{eqnarray*}}
\newcommand{\ean}{\end{eqnarray*}}

 \makeindex

\begin{document}

\title{Chambers's formula for  the graphene and the  Hou model with kagome periodicity and applications.} 
\author{Bernard Helffer ${}^{1,2}$, Philippe Kerdelhu\'e${}^{1}$ \& Jimena Royo-Letelier${}^{3}$ }
 \affil{ \footnotesize ${}^1$  D\'epartement de Math\'ematiques, CNRS UMR 8628, \\ F-91405 Orsay Cedex, France  \\ Bernard.Helffer@math.u-psud.fr, Philippe.Kerdelhue@math.u-psud.fr\\
 \footnotesize ${}^2$ Laboratoire de math\'ematiques J. Leray, univ. Nantes\\
\footnotesize ${}^3$  jimena.royo-letelier@m4x.org}
\date{}
\maketitle

 \begin{abstract}
 The aim of this article  is to prove that for the graphene model like for
 a model considered by the physicist Hou on a kagome lattice, there
exists a formula which is similar to the one obtained by Chambers  for
the Harper model.  As an application, we propose  a semi-classical
analysis of the spectrum of the Hou butterfly near a flat band.
\end{abstract}

  \section{Introduction}
 \subsection{A brief historics}
  Starting from the middle of the fifties \cite{Ha}, solid state physicists have been interested in the flux effects created by a magnetic field (see in the sixties Azbel \cite{Az}, Chambers \cite{Ch}) . In 1976  a  celebrated butterfly was proposed by D. Hofstadter \cite{Hof} to describe as a function of the flux $\gamma$  the spectrum (at the bottom) of a Schr\"odinger operator with constant magnetic field and periodic electric potential. About ten years later mathematicians start to propose rigorous proofs for this approximation and to analyze the model itself. The celebrated ten martinis conjecture about the Cantor structure when $ \gamma/2\pi$ is irrational was formulated by M. Kac and only solved a few years ago (see \cite{AJ} and references therein). We refer also to the survey of J. Bellissard \cite{bel} for a state of the art in 1991. Once a semi-classical (or tight-binding) approximation is done, involving a tunneling analysis
   we arrive (modulo a controlled smaller error) in the case of a square lattice to the so-called Harper model, which is defined on $\ell^2(\mathbb Z^2,\mathbb C)$ by
   $$
  ( H u)_{m,n} := \frac 12 (u_{m+1,n} + u_{m-1,n}) + \frac 12 e^{i  \gamma m} u_{m,n+1} +  \frac 12 e^{-i \gamma m} u_{m,n-1}\,,
  $$
  where $\gamma$ denotes the flux of the constant magnetic field through the fundamental cell of the lattice.\\
  When $\frac{\gamma}{2 \pi}$ is a rational,  a Floquet theory permits to show that the spectrum is the union   of the spectra of a family of  $q\times q$ matrices depending 
  on a parameter $\theta=(\theta_1,\theta_2) \in \mathbb R^2$.\\
  More precisely, when
\begin{equation}
\gamma=2\pi p/q\,,
\end{equation}
where $p\in\mathbb Z$ and $q\in\mathbb N^*$ are relatively prime, 
the two following matrices play an important role:
\begin{equation}
J_{p,q}={\rm diag}(e^{i(j-1)\gamma})\,,
\end{equation}
and
\begin{equation}
(K_q)_{jk}= 1\; \mbox{ if } k\equiv j+1\, [q]\,,\,0 \mbox{ else.}\,
\end{equation}
In the case of Harper, the family of matrices is 
\begin{equation}\label{bs1}
M_H(\theta_1,\theta_2) = 
\frac{1}{2} ( e^{i\theta_1} J_{p,q} + e^{-i\theta_1} J_{p,q}^* + e^{i\theta_2} K_q + e^{-i \theta_2} K_q^*)\,.
\end{equation}
The Hofstadter butterfly is then obtained as a picture in the rectangle $[-2,+2]\times [0,1]$ (see Figure 1). A point $(\lambda,\gamma/2\pi)$ is in the picture if there exists $\theta$ such that \break 
$\det (M_H(\theta_1,\theta_2)-\lambda) =0$ for some $\frac p q$ with $p/q \in [0,1]$ ($q\leq 50$).\\
The Chambers formula gives a very elegant formula for this determinant:
\begin{equation}
\det (M_H(\theta_1,\theta_2)-\lambda)  = f^H_{p,q}(\lambda)+ (-1)^q \left(\cos q\theta_1 + \cos q \theta_2\right)\,,
\end{equation}
where $f^H$ is a polynomial of degree $q$.\\

Many other models have been considered.  In the case of a triangular lattice, the second model is, according to \cite{Ke} (see also \cite{AKY}),
\begin{equation}\label{bs5}
M_T(\theta_1,\theta_2,\phi) = e^{i\theta_1} J_{p,q} + e^{-i\theta_1} J_{p,q}^* +  e^{i\theta_2} K_q + e^{-i \theta_2} K_q^* +   e^{i\phi} e^{i(\theta_1-\theta_2)} J_{p,q}K_q^*  +
 e^{-i \phi} e^{i(\theta_2-\theta_1)} K_q  J_{p,q}^*
\end{equation}
with $ \phi=-\gamma/2$.\\
The Chambers formula in this case takes the form
\begin{equation}
\det (M_T(\theta_1,\theta_2,\phi)-\lambda)  = f^T_{p,q,\phi}(\lambda)+ (-1)^{q+1} \left(\cos q\theta_1 + \cos q \theta_2 + \cos q (\theta_2-\theta_1 -\phi )\right)\,.
\end{equation}
The resulting spectrum is given in Figure 2.\\

In the case of the hexagonal lattice, which appears also in the analysis of the graphene, we have to analyze
\begin{equation}\label{bsG}
M_G (\theta_1,\theta_2):= \left(
\begin{array}{cc}
0& I_q + e^{i\theta_1} J_{p,q}+ e^{i \theta_2} K_q\\
I_q + e^{-i\theta_1} J_{p,q}^* + e^{-i \theta_2} K_q^* & 0\end{array}
\right)
\end{equation}
We denote by $P_G$ the characteristic polynomial of $M_G$.
The resulting spectrum is given in Figure 3.\\
Finally, inspired by the physicist Hou, P. Kerdelhu\'e and J. Royo-Letelier \cite{KR} have shown that for the kagome lattice,  the following approximating model is relevant: 
we consider the matrix
\begin{equation}\label{1.9}
 M_K(\theta_1,\theta_2,\omega)= \left(\begin{array}{ccc}
0 & A(\theta_1,\theta_2,\omega) & B(\theta_1,\theta_2,\omega)\\
A^*(\theta_1,\theta_2,\omega) & 0 & C(\theta_1,\theta_2,\omega)\\
B^*(\theta_1,\theta_2,\omega) & C^*(\theta_1,\theta_2,\omega) & 0
\end{array}\right)\,,
\end{equation}
with 
\ban
A(\theta_1,\theta_2,\omega)&=&e^{i(\omega+\frac\gamma{8})} (e^{-i\theta_1}J_{p,q}^*+e^{-i\frac\gamma2 }e^{-i(\theta_1-\theta_2)}J_{p,q}^*K_q)\\
B(\theta_1,\theta_2,\omega)&=&e^{-i(\omega+\frac{\gamma}{8})} (e^{-i\theta_1}J_{p,q}^*+e^{-i\theta_2}K_q^*)\\
C(\theta_1,\theta_2,\omega)&=&e^{i(\omega+\frac{\gamma}{8})} (e^{-i\frac\gamma2 }e^{i(\theta_1-\theta_2)}J_{p,q}K_q^*+e^{-i\theta_2}K_q^*)\,.\ean
 Here $\omega$ is a parameter appearing in the model (most of the physicists consider without justification the case $\omega=0$). We refer to \cite{KR} for a discussion of this point.\\

The trigonometric polynomial
\begin{equation}
 (x,\xi) \mapsto p^\triangle(x,\xi)=\cos x + \cos \xi +\cos (x-\xi)
 \end{equation}
which was playing an important role in the analysis of the triangular Harper model  (see Claro-Wannier \cite{CW}  and  Kerdelhu\'e \cite{Ke}) will also appear in our analysis.\\
We denote by $P_K(\theta_1,\theta_2,\omega,\lambda)$ the characteristic polynomial
$\det(M_K(\theta_1,\theta_2,\omega)_\lambda)$.\\
 
 \subsection{Main results}
 The aim of this article is to prove that, for a model considered by Hou \cite{Hou}, there exists a formula which is similar to the one obtained by Chambers \cite{Ch} for the Harper model. (see also Helffer-Sj\"ostrand \cite{HS1}, \cite{HS2}, Bellissard-Simon \cite{belsim}, C. Kreft \cite{Kr}, I. Avron (and coauthors) \cite{AKY}). Such an existence was motivated by computations of \cite{KR}. We also consider the case of the graphene, where a huge litterature in Physics exists  (see \cite{DM} and references therein) which is sometimes unaware of semi-classical mathematical results of the nineties. Note that the Chambers formula plays an important role in the semi-classical analysis of the Harper's model (see for example \cite{HS2}). \\

The first statement is probably well known in the physical literature.
  \begin{theorem}[Graphene]\label{th1}
\begin{equation}\label{aprioriG}
P_G(\theta_1,\theta_2,\lambda)= (-1)^q\,  {\rm det} (M_T(\theta_1,\theta_2,0) +3  -\lambda^2) \,.
\end{equation} 
\end{theorem}
The second statement was to our knowledge unobserved.
\begin{theorem}[Kagome]\label{th2}~\\
For any $\omega$, there exists a polynomial $Q_\omega$ of degree $3q$, with real coefficients, depending on $p,q$, such that
\begin{equation}\label{aprioriK}
P_K(\theta_1,\theta_2,\omega,\lambda)=Q_\omega(\lambda)+2 p^\triangle(q(\theta_1+p\pi),q(\theta_2+p\pi)) R_\omega(\lambda)\,,
\end{equation}
with
\begin{equation}\label{defRommega}
R_\omega(\lambda):=\left(\lambda+2\cos(3\omega-\frac{\gamma}{8})\right)^q\,.
\end{equation} 
Moreover the principal term of $Q_\omega(\lambda)$ is $\lambda^{3q}$.
\end{theorem}
We call $k$-th band the set described when $(\theta_1,\theta_2) \in \mathbb R^2$ by the $k$-th eigenvalue of the matrix $M_K$. We will call this band flat if this $k$-th eigenvalue is independent of $(\theta_1,\theta_2)$.
\begin{corollary}
A flat band exists if and only if
$$
 Q_\omega(-2\cos(3\omega-\frac{\gamma}{8}))=0\,.
$$
\end{corollary}
\begin{remark}~
\begin{itemize}
\item 
$Q_\omega$ is a trigonometric potential in $3\omega$.
\item For $(p,q)$ given, the set of the $\omega$'s such that a flat band exists is discrete.
 Formula (\ref{1.9}) shows indeed that the expression 
$P_K(\theta_1,\theta_2,\omega,-2\cos(3\omega-\gamma/8))$, which according to Theorem \ref{th2} is independent of  $(\theta_1,\theta_2)$, takes the form $\Sigma_{j=-9q}^{9q} a_j e^{i j \omega}$  with $a_{9q}=e^{-i\frac{3 \gamma q}{8}}$. 
\end{itemize}
\end{remark}
\subsection{Examples}
Let us illustrate by some examples mainly extracted of \cite{KR}.\\
In the case when $q=1$ and $p=0$, one finds, for the Hou's model:
$$
P_K(\theta_1,\theta_2,\omega,\lambda)= \lambda^3- 6  \lambda   -  4 \cos(3\omega)  - 2  \left(   \lambda+ 2 \cos(3\omega) \right) p^\triangle (\theta_1,\theta_2) \,.
$$
Hence, we have in this case:
$$
Q_\omega(\lambda)= \lambda^3- 6  \lambda   - 4 \cos(3 \omega) \,.
$$
It is then natural to ask if the two polynomial have a common zero. The condition reads:
$$
Q_\omega(- 2 \cos(3 \omega) \ ) =0\,.
$$
We get:
$$
 (\cos 3 \omega)^3  - \cos 3 \omega =0\,,
$$
hence $\cos 3 \omega =0$ or $\cos 3 \omega = \pm 1$. So a  "flat band" appears when  $\omega =0$, which was mostly considered in the physical literature. Note that in \cite{KR}, it is proved only that $\omega\rightarrow 0$ as a function of the initial semi-classical parameter.
The set of $\omega$'s for which we have a flat band is $\{\omega_k =k\frac \pi 6\,,\, k\in \mathbb Z\}\,.$ 
\\

Another example is, as shown in \cite{KR} (Proposition 1.13), for $\omega=\pi/8$ and $p/q=3/2$. The bands are $\{-2\}$ (with multiplicity $2$), $[1-\sqrt{6},1-\sqrt{3}]$, $[1-\sqrt{3},1]$, $[1,1+\sqrt{3}]$ and 
$[1+\sqrt{3},1+\sqrt{6}]$.

\subsection{Organization of the paper}
This paper is organized as follows. In Section \ref{ssymm} we establish symmetry properties of the two matrices $J_{p,q}$ and $K_q$. In Section \ref{sBeSi} we recall
 how a method due to Bellissard-Simon permits to establish the Chambers formula for a square lattice or a triangular lattice. In Section \ref{sgraphene}, we give an application to the case of the graphene.
 Section \ref{sHou} is devoted to the proof of the main theorem for the kagome lattice.
 In Section \ref{soverlap}, we establish the non overlapping of the bands in the case of the kagome lattice. Section \ref{s7} gives as an application a semi-classical analysis near a flat band and we finish with a conclusion.

\section{Symmetries}\label{ssymm}
We recall some basic symmetry properties of the two matrices $J_{p,q}$ and $K_q$. Some of them were used in the previous literature, some other are new.
We first recall that
\begin{equation}\label{comm1}
J_{p,q } K_q = \exp ( -2i \pi \frac p q )\, K_q J_{p,q} \,.
\end{equation}
and (take the complex conjugation and the adjoint )
\begin{equation}\label{comm2}
K_q^* J_{p,q} = \exp ( - 2i \pi \frac p q )\, J_{p,q} K_q^* \,.
\end{equation}

\begin{lemma}\label{lemme1.2}
There exist unitary matrices   $U$ and  $V$ in $M_q(\mathbb C$) such that
\ba \label{f2}  U^*K_q^*U&=&J_{p,q}  \\
\label{f3}  U^*J_{p,q}U&=&K_q \\
\label{f4} V^*K_q^*V&=&J_{p,q}\\
\label{f5}  V^*J_{p,q}V&=&(-1)^p e^{-i\frac{\gamma}{2}}J_{p,q}^*K_q\\
\label{f6} V^*((-1)^p e^{-i\frac{\gamma}{2}}J_{p,q}^*K_q)V&=& K_q^*  \ea
\end{lemma}
\begin{remark} Note from (\ref{comm1}) and (\ref{comm2}) that the pairs $(J_{p,q},K_q)$ and $(K_q^*,J_{p,q})$ satisfy the same commutation relation. 
(\ref{f2}) et (\ref{f3})  make explicit the unitary equivalence between this representation and the one used in \cite{KR}.
\end{remark}

{\bf Proof}\\
$U$ is actually the discrete Fourier transform:
 \begin{equation}
 U_{j,k}=q^{-1/2} e^{-i\gamma(j-1)(k-1)}\,,\, j,k=1,\cdots, q\,.
 \end{equation}
It is easy to verify (\ref{f2}) et (\ref{f3}).\\
For (\ref{f4}), we observe that, $J_{p,q}$  being diagonal,  (\ref{f4}) is verified for any matrix $V$ in the form
 $$ V=UD \,,$$
where $D$ is a diagonal unitary matrix
$$
D={\rm diag}(d_j)\,,
$$
with $|d_j|=1$.\\
We are looking for the $d_j$'s and a complex number   $c$ of module $1$ such that
$$V^*J_{p,q}V=cJ_{p,q}^*K_q\,.$$

If we think of the indices as elements in $\mathbb Z/q \mathbb Z$, we have:
$$(V^*J_{p,q}V)_{j,k}=d_{j+1}\bar d_j \delta_{j+1,k}\,,$$
and
$$(J_{p,q}^*K_q)_{j,k}=e^{-i(j-1)\gamma} \delta_{j+1,k}\,.$$
We want to have 
 $$d_1=1\,, \,d_ {j+1}=c\,e^{-i(j-1)\gamma} \mbox{ for } j>1\,,$$ 
but also: $$d_{q+1}=1\,.$$
This implies
$$e^{-i\gamma\frac{q(q-1)}{2}}c^q=1\,.
$$
So we choose
$$c=e^{i\gamma\frac{q-1}{2}}=(-1)^p e^{-i\frac{\gamma}{2}}\,.
$$
We then obtain
 $$ V^*(J_{p,q}^*K_q)V=\bar c K_q^*J_{p,q}J_{p,q}^*=\bar c K_q^*\,.
$$
\qed
\\

\section{Harper on square and triangular lattice}\label{sBeSi}
 We recall in this section  the approach of Bellissard-Simon \cite{belsim}, initially introduced for the analysis of the Harper model, we apply it for the case of the triangular lattice. Note that this second situation was recently analyzed in \cite{AKY} and \cite{AEG}.
\subsection{The case of Harper}

We start from the general formula
\begin{equation} \label{bs2}
\det (M -\lambda I_q) = (-\lambda)^q \exp \Tr \left(\log (I_q- \frac{M}{\lambda})\right)\,.
\end{equation}
This implies
\begin{equation} \label{bs3} 
\det (M_H(\theta_1,\theta_2) -\lambda I_q) = (-\lambda)^q \exp \left( - \sum_{k\geq 1}  \lambda^{-k}\,  \frac{\Tr M_H(\theta_1,\theta_2)^k}{k } \right)\,.
\end{equation}
The next point is to observe that
\begin{equation}\label{bs4}
\Tr (J_{p,q}^{\ell_1} K_q^{\ell_2}) =0\,,\, \mbox{ except } \ell_1 \equiv 0 \mbox{ and } \ell_2 \equiv 0 \mbox{ mod q}\,.
\end{equation}
The only term which depends on $(\theta_1, \theta_2)$ in 
 $\frac{1}{k \, \lambda^k}\, \Tr M_H(\theta_1,\theta_2)^k$ (for $k\leq q$) corresponds to $k=q$ and is simply: $ \frac{2}{\lambda^q}\,\left( \cos q \theta_1 + \cos q \theta_2\right)$.\\
 
 The general term is indeed 
 $$\exp i \left(\ell_1 \theta_1 - \ell_1^* \theta_1 + \ell_2 \theta_2-\ell_2^* \theta_2\right) \;  \Tr J_{p,q}^{\ell_1-\ell_1^*} K_q^{\ell_2-\ell_2^*}\,,
 $$
 with $\ell_1\geq 0\,,$   $\ell_1^*\geq 0\,, $ $\ell_2\geq 0\,,$ $\ell_2^*\geq 0\,,$  and  $\ell_1 + \ell_1^* + \ell_2 + \ell_2^* \leq q\,$.\\
 But \eqref{bs4} implies that the non vanishing terms (depending effectively on $(\theta_1,\theta_2)$) can only correspond to 
 $$
 \ell_1\equiv \ell_1^*\,\mbox{ and } \ell_2 \equiv \ell_2^*\,, \mbox{ with } |\ell_1 -\ell_1^*| + |\ell_2 - \ell_2^*|\neq 0\,.
 $$
 A case by case analysis leads to only four non zero terms corresponding to \break $\ell_1=q, \ell^*_1=0, \ell_2=0, \ell_2^*=0$,  and the three permutations of this case.
 Hence we have proved:
 \begin{proposition}
 \begin{equation}
 \det (M_H(\theta_1,\theta_2) -\lambda I_q)  = f^H_{p,q} (\lambda) + (-1)^{q+1} 2 \left(\cos q \theta_1 + \cos q \theta_2\right)\,.
 \end{equation}
 \end{proposition}
\subsection{The case of  Harper on a triangular lattice} 
We first treat the case with $\phi$ as a free parameter.\\
The starting point is the same but this time  the general term is
 $$
 \exp i \left(\ell_1 \theta_1 - \ell_1^* \theta_1 + \ell_2 \theta_2-\ell_2^* \theta_2 + (\ell_3-\ell_3^*)  (\theta_1-\theta_2) \right) \;  \Tr J_{p,q}^{\ell_1-\ell_1^*+ \ell_3-\ell_3^*} K_q^{\ell_2-\ell_2^*-\ell_3 + \ell_3^*} \,,
 $$
  with $\ell_1\geq 0\,,$   $\ell_1^*\geq 0\,,$ $\ell_2\geq 0\,,$ $\ell_2^*\geq 0\,,$ $\ell_3\geq 0\,,$ $\ell_3^*\geq 0\,,$ and  
  \begin{equation}\label{cut}
  \ell_1 + \ell_1^* + \ell_2 + \ell_2^*+ \ell_3 + \ell_3^* \leq q\,.
  \end{equation}
But \eqref{bs4} implies that the non vanishing terms can only correspond to 
 \begin{equation}\label{cut2}
 \ell_1-\ell_1^*+ \ell_3-\ell_3^* \equiv 0\,\mbox{ and } \ell_2-\ell_2^*-\ell_3 + \ell_3^* \equiv  0 \,,
 \end{equation}
 with
 \begin{equation}\label{cut3}
  \ell_1-\ell_1^*+ \ell_3-\ell_3^* \neq 0\,\mbox{ or } \ell_2-\ell_2^*-\ell_3 + \ell_3^* \neq  0 \,.
 \end{equation}
 We have six evident cases corresponding to all indices equal to $0$ except one equal to $q$. It remains  to discuss if   there are  other cases.\\
 We introduce the auxiliary parameters:
 $$
 \tilde \ell_1 = \ell_1 + \ell_3\,,\, \tilde \ell_1^* = \ell_1^* + \ell_3^*\,,\,  \tilde \ell_2 = \ell_1 + \ell_3^*\,,\, \tilde \ell_1^* = \ell_2^* + \ell_3\,,
 $$
and with these conditions we get:
\begin{equation}\label{cut2prime}
\tilde  \ell_1 -\tilde \ell_1^* \equiv 0\,\mbox{ and } \tilde \ell_2- \tilde \ell_2^* \equiv  0 \,,
 \end{equation}
 with
 \begin{equation}
 \tilde  \ell_1 -\tilde \ell_1^* \neq 0\,\mbox{ or } \tilde \ell_2- \tilde \ell_2^* \neq  0 
 \end{equation}
 This looks rather similar to the previous situation except the bounds on the $\tilde \ell_j$.\\
 In the case by case discussion, we first verify that for each congruence it is enough (using \eqref{cut})  to look at $\tilde \ell_j -\tilde \ell_j^*= -q,0,q$ hence to nine cases  but the second condition eliminates one case.  One can also eliminate two cases corresponding to $(\tilde \ell_1 -\tilde \ell_1^*)(\tilde \ell_2 -\tilde \ell_2^*)>0$ using again the condition \eqref{cut}. Hence it remains six cases, each one containing  one of the evident cases.\\ 
 Let us look at one of these six cases: 
 $$\tilde\ell_1 = \tilde \ell_1^* +q \,,\, \tilde \ell_2 = \tilde \ell_2^* -q\,.
 $$
This reads
$$
\ell_1 + \ell_3 = \ell_1^* + \ell_3^* +q \,,\, \ell_2+\ell_3^* = \ell_2^* + \ell_3 -q\,.
$$
The left part together with \eqref{cut} implies $\ell_1^*=\ell_3^*=0$ and the right part implies $\ell_2=0$. Hence it remains:
$$
\ell_1 + \ell_3 = q\,,\, \ell_2^* = q-\ell_3 = \ell_1\,.
$$ 
Using again the condition on the sum we get  $\ell_2^* =\ell_1 =0\,$, hence finally $\ell_3=0\,$. We are actually in one of the six announced  trivial cases.\\

 \begin{proposition}
 \begin{equation}
 \det (M_T(\theta_1,\theta_2, \phi ) -\lambda I_q)  = f^{T}_{p,q,\phi} (\lambda) + (-1)^{q+1} 2 \left(\cos q \theta_1 + \cos q \theta_2 + (-1)^{q+1} \cos q (\theta_1-\theta_2+ \phi)\right)\,.
 \end{equation}
 \end{proposition}
 What remains is to compute the coefficients in the six cases (actually three cases are enough because the sum should be real). We only compute the new case.
 As  $$((-1)^p e^{-i\gamma/2}J_{p,q}K_q^*)^q=I_q$$ we immediately get as coefficient
$\cos(q\theta_1)+\cos(q\theta_2)+(-1)^{pq}\cos(q\theta_1-q\theta_2+\pi p+q\phi)$ which can be written observing that 
 $(-1)^{(p+1)(q+1)}=1$ ($p$ and  $q$ being mutually prime):
$$\cos(q\theta_1)+\cos(q\theta_2)+(-1)^{q+1}\cos(q\theta_1-q\theta_2+q\phi)\,.$$

\begin{remark}
Similar formulas appear in \cite{AEG}.
\end{remark}

 \section{The hexagonal or graphene case}\label{sgraphene}

Taking the square of the matrix given by (\ref{bsG}), we obtain
\begin{equation}\left(
\begin{array}{cc} 3 I_q + M_T(\theta_1,\theta_2,0) 
 & 0\\ 0 &   3 I_q +  \hat M_T(\theta_1,\theta_2,0) 
\end{array}
\right)
\end{equation}
with
\begin{multline}
\hat M_T(\theta_1,\theta_2,0) =  e^{i\theta_1} J_{p,q} +  e^{-i\theta_1} J_{p,q}^* + e^{i \theta_2} K_q + e^{-i \theta_2} K_q^*\\+  e^{i(\theta_1-\theta_2)} K_q^* J_{p,q} + e^{-i(\theta_1-\theta_2)} J_{p,q}^* K_q \,.
\end{multline}

For the second term  we have just an exchange of $J_{p,q}$ and $K_q$. It is clear by supersymmetry that the two terms have the same non-zero eigenvalues. If we control  the multiplicity 
 this will give the isospectrality. If we  introduce 
 $$
 \mathcal A =   I_q + e^{i\theta_1} J_{p,q}+ e^{i \theta_2} K_q \,,
 $$
 the two operators read $\mathcal A \mathcal A^*$ and $\mathcal A^* \mathcal A\,$. \\
 Consider indeed $u\neq 0$ such that
 $$
 \mathcal A \mathcal A^* u = \lambda u\,.
 $$
 Then we get
 $$
 \mathcal A^* \mathcal A \mathcal A^* u = \lambda \mathcal A^* u\,.
 $$
 If $\lambda \neq 0$, then $\mathcal A^*  u \neq 0$ and is consequently an eigenvector of  $ \mathcal A^* \mathcal A $. The multiplicity is also easy to follow.

 Hence we get easily an equation for the square of the eigenvalues. But it has been shown in \cite{KR} (by conjugation by $\left(\begin{array}{cc} -I_q&0\\0&I_q \end{array}\right) $ ), 
  that the spectrum is invariant by  $\lambda \rightarrow -\lambda$. Hence looking at the first characteristic polynomial  gives us all the squares of the eigenvalues of $M_G + 3 I_{2q}$, counted with multiplicity.

So we have proved Theorem \ref{th1}. 
Hence the spectrum will consists of $q$ bands in $\mathbb R^+$ and of $q$ bands in $\mathbb R^-$ obtained by symmetry. We will show in the next section that these bands are not overlapping but that possibly touching. 
The last (maybe standard) observation is  that  the two central gaps for the Graphene-model  are effectively touching at $0$.
 We have to show that $0$ belongs to the spectrum :
 \begin{proposition}
  There exists $(\theta_1,\theta_2)\in\mathbb R^2$ such that
$$
{\rm det} (M_G(\theta_1,\theta_2))  =0\,.
$$
\end{proposition}
It is actually enough to show:
\begin{lemma} 
 There exists $(\theta_1,\theta_2)\in\mathbb R^2$ such that
$$\det(I_q+e^{i\theta_1} J_{p,q}+e^{i\theta_2} K_q)=0\,.$$
\end{lemma}
{\bf Proof}\\
We consider the polynomial
\begin{multline}
P(\lambda)=\det(-\lambda\,I_q+e^{i\theta_1} J_{p,q}+e^{i\theta_2} K_q)=\\
\det
\left(
\begin{array}{ccccc}
 -\lambda+e^{i\theta_1} & e^{i\theta_2} & 0 & \cdots & 0 \\
0 & -\lambda+e^{i2\pi p/q}e^{i\theta_1} & e^{i\theta_2} & \cdots & 0 \\
\vdots & \vdots & \vdots & \ddots & \vdots \\
e^{i\theta_2} & 0 & 0 & \cdots & -\lambda+e^{i2\pi p (q-1)/q}e^{i\theta_1}                                                                    
      \end{array}\right )\end{multline}

$P$ has degree $q$, the coefficient of $\lambda^q$ is $(-1)^q$, and $$P(\lambda)=(-1)^{q-1}e^{iq\theta_2}$$ if
$\lambda=e^{i2\pi k/q}e^{i\theta_1}$ for $k\in\{0,\cdots,q-1\}$, i.e. if $\lambda^q=e^{i q\theta_1}$. Hence
$$P(\lambda)=(-1)^q(\lambda^q-e^{i q\theta_1}-e^{iq\theta_2})\,.$$
Considering $\lambda=-1$ gives $$\det(I_q+e^{i\theta_1} J_{p,q}+e^{i\theta_2} K_q)=1-e^{i q(\theta_1+\pi)}-e^{iq(\theta_2+\pi)}\,.$$ The choice of 
$\theta_1 =\pi+\pi/(3q)$ and $\theta_2=\pi-\pi/(3q)$ achieves the proof. \qed

\begin{remark}
Interesting new results concerning the graphene case and the computation of Chern classes have been obtained recently in \cite{AEG} and \cite{AKY}.
\end{remark}

\section{Proof of Theorem \ref{th2}}\label{sHou}
Although the Bellissard-Simon approach gives a partial proof of Theorem \ref{th2}, the proof given below goes much further by implementing the symmetry considerations described in Section \ref{ssymm}.
\subsection{First a priori form}
We first establish:
\begin{lemma}\label{lem1}
There exist polynomials $T_{j,k}$, $-1\leq j,k\leq1$ such that, for all $(\theta_1,\theta_2)\in \mathbb R^2$ 
\be \label{f1} P_K(\theta_1,\theta_2,\omega,\lambda)=\sum_{j,k\in\{-1,0,1\}} e^{i(q(j\theta_1+k\theta_2))} T_{j,k}(\lambda)\,. \ee
\end{lemma}

{\bf Proof:}\\
 We define the matrix $S(\theta_1,\theta_2)$, which is  unitary equivalent with $M_K(\theta_1,\theta_2,\omega)$,  by
\begin{multline} \label{defS} 
S(\theta_1,\theta_2)=
\left(\begin{array}{ccc}
         e^{-i\theta_1}J_{p,q}^* & 0 &0\\
0 & I_q & 0\\
0 & 0 & e^{i\theta_2}K_q
        \end{array}\right)^*
M_K(\theta_1,\theta_2,\omega)\left(\begin{array}{ccc}
         e^{-i\theta_1}J_{p,q}^* & 0 &0\\
0 & I_q & 0\\
0 & 0 & e^{i\theta_2}K_q
        \end{array}\right)\,.
        \end{multline}
        A computation shows that
 \begin{multline} \label{proprS} 
S(\theta_1,\theta_2)=
       \left(\begin{array}{ccc}
    0 &    e^{i(\omega+\frac{\gamma}{8})} (I_q+e^{-i\frac\gamma2 }e^{i\theta_2}K_q)   & e^{-i(\omega+\frac{\gamma}{8})} (e^{i\theta_2}K_q+e^{i\theta_1}J_{p,q})\\
e^{-i(\omega+\frac{\gamma}{8})} (I_q+e^{i\frac\gamma2 }e^{-i\theta_2}K_q^*) & 0 & e^{i(\omega+\frac{\gamma}{8})} (e^{-i\frac\gamma2 }e^{i\theta_1}J_{p,q}+I_q)\\
  e^{i(\omega+\frac{\gamma}{8})} (e^{-i\theta_2}K_q^*+e^{-i\theta_1}J_{p,q}^* ) & e^{-i(\omega+\frac{\gamma}{8})} (e^{i\frac\gamma2 }e^{-i\theta_1}J_{p,q}^*+I_q)  & 0   \end{array}\right)\,.
\end{multline}
Hence $M_K(\theta_1,\theta_2,\omega)$ and $S(\theta_1,\theta_2)$ have the same characteristic polynomial and  coming back to the definition of the
determinant, we can verify  that $P$ is a polynomial of degree  $q$ in $(e^{-i\theta_1},e^{i\theta_1})$, and also of  degree  $q$ in
$(e^{-i\theta_2},e^{i\theta_2})$.\\
Then we observe that
$$\left(\begin{array}{ccc}
         J_{p,q}&0&0\\
0&J_{p,q}&0\\
0&0&J_{p,q}
        \end{array}\right)^*
M_K(\theta_1,\theta_2,\omega)
\left(\begin{array}{ccc}
         J_{p,q}&0&0\\
0&J_{p,q}&0\\
0&0&J_{p,q}
        \end{array}\right)=M_K(\theta_1,\theta_2+\frac{2\pi p}{q},\omega)$$
and 
$$\left(\begin{array}{ccc}
         K_q&0&0\\
0&K_q&0\\
0&0&K_q
        \end{array}\right)^*
M_K(\theta_1,\theta_2,\omega)
\left(\begin{array}{ccc}
         K_q&0&0\\
0&K_q&0\\
0&0&K_q
        \end{array}\right)= M_K(\theta_1-\frac{2\pi p}{q},\theta_2,\omega)\,.$$
        As $P_K$ is $2\pi$-periodical in $\theta_1$ and  $\theta_2$, and  $p$ et $q$ are mutually prime, $P$ is\footnote{This argument is already present in a similar context in \cite{belsim}.}
$(2\pi/q)$-p\'eriodical in $\theta_1$ and $\theta_2$. One can indeed use B\'ezout's theorem observing that $1 = u p + v q$ (with  $u$ and  $v$ in $\mathbb Z$), hence  $\frac 1 q = u\,  \frac p q + v$.\\
\\ \qed\\

\subsection{Improved a priori form}~\\
Here we prove the existence of two polynomials $Q_\omega$ and $R_\omega$, with real coefficients,
depending on $\gamma$ and possibly on $\omega$, but not on
$(\theta_1,\theta_2,\omega)$, such that
\begin{equation}\label{apriori}
P_K(\theta_1,\theta_2,\omega,\lambda)=Q_\omega(\lambda)+p^\triangle(q(\theta_1+p\pi),q(\theta_2+p\pi)) R_\omega(\lambda)\,.
\end{equation}

In view of Lemma \ref{lem1}, it remains to prove that $P(\theta_1+p\pi,\theta_2+p\pi)$ is invariant by the "rotation  of angle  $-2\pi/3$" $r$ which leaves invariant  $p^\triangle$
and is defined by $$r(\theta_1,\theta_2)=(-\theta_1+\theta_2,-\theta_1)\,,
$$
 and by the symmetry   $s$ defined by
 $$ s(\theta_1,\theta_2)=(\theta_2,\theta_1)\,.
 $$
We now introduce
 \begin{equation}
 N(\theta_1,\theta_2)=(-1)^p M_K(\theta_1+p\pi,\theta_2+p\pi,\omega)\,,
 \end{equation}
 and 
 \begin{equation}
 L_{p,q}=(-1)^pe^{-i\frac{\gamma}{2}}J_{p,q}^*K_q\,.
 \end{equation}
With this notation and $\omega'=\omega+\gamma/8$, $N(\theta_1,\theta_2)$ reads:
\begin{equation}
\left(\begin{array}{ccc}
0 &    e^{i\omega'} (e^{-i\theta_1}J_{p,q}^*+e^{-i(\theta_1-\theta_2)}L_{p,q}) & e^{-i\omega'} (e^{-i\theta_1}J_{p,q}^*+e^{-i\theta_2}K_q^*)  \\
 e^{-i\omega'} (e^{i\theta_1}J_{p,q}+e^{i(\theta_1-\theta_2)}L_{p,q}^*) &  0 & e^{i\omega'} (e^{i(\theta_1-\theta_2)}L_{p,q}^*+e^{-i\theta_2}K_q^*)\\
e^{i\omega'} (e^{i\theta_1}J_{p,q}+e^{i\theta_2}K_q) & e^{-i\omega'} (e^{-i(\theta_1-\theta_2)}L_{p,q}+e^{i\theta_2}K_q)  & 0
               \end{array}
\right)
\end{equation}

We will show that the characteristic polynomial of $N$ is invariant by $r$ and $s$. We have seen that
$$V^*K_q^*V=J_{p,q}\,, \, V^*J_{p,q}V=L_{p,q} \mbox{ and  } V^*L_{p,q}V=K_q^*\,.$$
We easily see that~:\\
\begin{lemma}
\begin{equation}
\left(\begin{array}{ccc}
         0&V&0\\
0&0&V\\V&0&0
        \end{array}\right)^* N(r(\theta_1,\theta_2)) \left(\begin{array}{ccc}
         0&V&0\\
0&0&V\\V&0&0
        \end{array}\right)=N(\theta_1,\theta_2)\,.
        \end{equation}
  \end{lemma}
Hence the characteristic polynomial is invariant by $r$.\\
We have already used that $\bar K_q=K_q$ et $\bar J_{p,q}=J_{p,q}^*$ and we have consequently~:
\ben  U^*\overline{J_{p,q}} U=K_q^*,~~U^*\overline{K_q} U=J_{p,q}^* \mbox{ and  }  U^*\overline{L_{p,q}} U=L_{p,q} \een
It is then easy to get:\\
\begin{lemma}
 $$\left(\begin{array}{ccc}
         0&0&U\\
0&U&0\\
U&0&0
\end{array}\right)^*
\overline{N(\theta_2,\theta_1)}
\left(\begin{array}{ccc}
         0&0&U\\
0&U&0\\
U&0&0
\end{array}\right)
=N(\theta_1,\theta_2)\,.$$
\end{lemma}
Hence the characteristic polynomial is invariant by $s$.\\
\qed

\subsection{End of the proof}~\\
We now make explicit the polynomial $R_\omega$. (\ref{defS}) reads:
  \begin{multline}
 2\,e^{iq(\theta_1-\theta_2)} (Q_\omega (\lambda)+(\cos(q(\theta_1-\theta_2))+(-1)^{pq}\cos(q\theta_1)+(-1)^{pq}\cos(q\theta_2))R_\omega(\lambda))\\
=2\det
\left(\begin{array}{ccc}
    -e^{-i\theta_2}\lambda\,I_q &    e^{i(\omega+\frac{\gamma}{8})} (e^{-i\theta_2} I_q+e^{-i\frac\gamma2 }K_q)   & e^{-i(\omega+\frac{\gamma}{8})} (K_q+e^{i(\theta_1-\theta_2)}J_{p,q})\\
e^{-i(\omega+\frac{\gamma}{8})} (I_q+e^{i\frac\gamma2 }e^{-i\theta_2}K_q^*) & -\lambda\,I_q & e^{i(\omega+\frac{\gamma}{8})} (e^{-i\frac\gamma2 }e^{i\theta_1}J_{p,q}+I_q)\\
  e^{i(\omega+\frac{\gamma}{8})} (e^{i(\theta_1-\theta_2)}K_q^*+J_{p,q}^* ) & e^{-i(\omega+\frac{\gamma}{8})} (e^{i\frac\gamma2 }J_{p,q}^*+e^{i\theta_1}I_q)  & -e^{i\theta_1}\lambda\,I_q
   \end{array}\right) \label{g1}
\end{multline} 
This equality between holomorphic functions holds for real $(\theta_1,\theta_2)$ and hence for complex $(\theta_1,\theta_2)$. Let $t$ be a real parameter and
take $\theta_1=-\theta_2=it$ in (\ref{g1}). The limit $t\rightarrow +\infty$ gives:
\ban R_\omega (\lambda)&=&2\,\det \left(\begin{array}{ccc}
   0 &    e^{i(\omega+\frac{\gamma}{8})} e^{-i\frac\gamma2 } K_q   &   e^{-i(\omega+\frac{\gamma}{8})}  K_q\\
e^{-i(\omega+\frac{\gamma}{8})} I_q & -\lambda\,I_q & e^{i(\omega+\frac{\gamma}{8})} I_q\\
  e^{i(\omega+\frac{\gamma}{8})} J_{p,q}^*  & e^{-i(\omega+\frac{\gamma}{8})} e^{i\frac\gamma2 } J_{p,q}^*  & 0
   \end{array}\right)\\
&=& 2\det\left(\begin{array}{ccc}
                K_q & 0 & 0\\
0 & I_q & 0\\
0 & 0 & J_{p,q}^*
               \end{array}\right)
\det \left(\begin{array}{ccc}
   0 &    e^{i(\omega+\frac{\gamma}{8})}e^{-i\frac\gamma2 } I_q   &   e^{-i(\omega+\frac{\gamma}{8})}  I_q\\
e^{-i(\omega+\frac{\gamma}{8})} I_q & -\lambda\,I_q & e^{i(\omega+\frac{\gamma}{8})} I_q\\
  e^{i(\omega+\frac{\gamma}{8})} I_q  & e^{-i(\omega+\frac{\gamma}{8})} e^{i\frac\gamma2 } I_q  & 0
   \end{array}\right) \,.\ean
$J_{p,q}$ and $K_q$ are conjugate, hence
\be  \det\left(\begin{array}{ccc}
                K_q & 0 & 0\\
0 & I_q & 0\\
0 & 0 & J_{p,q}^*
               \end{array}\right)=1 \,, \ee
and a straightforward computation gives
\be \det \left(\begin{array}{ccc}
   0 &    e^{i(\omega+\frac{\gamma}{8})} e^{-i\frac\gamma2 } I_q   &  e^{-i(\omega+\frac{\gamma}{8})} I_q\\
e^{-i(\omega+\frac{\gamma}{8})} I_q & -\lambda\,I_q & e^{i(\omega+\frac{\gamma}{8})} I_q\\
  e^{i(\omega+\frac{\gamma}{8})} I_q  & e^{-i(\omega+\frac{\gamma}{8})} e^{i\frac\gamma2 } I_q  & 0
   \end{array}\right) 
=\left(\lambda+2\cos(3\omega-\frac\gamma8)\right)^q\,. \ee \\ \qed

\section{On the non-overlapping of the bands}\label{soverlap}
The non overlapping of the bands has been proved in \cite{belsim} who refers for one part to a general  argument to Reed-Simon \cite{ReSi}. The fact that except at the center for $q$ even, the bands do not touch has been proven by P. Van Mouche \cite{VM}.  We show below that the non overlapping of the bands is a general property for all the considered domains but that the "non touching" property was specific of the Harper model.

 \begin{lemma}
Let $f(\lambda)$ be a real polynomial of degree $q$, such that, for any $\mu \in I=]a,b[$, $f(\lambda)=\mu$  has $q$ real solutions.  Then $f'(\lambda)\neq 0$, for any $\lambda$ such that $f(\lambda)=\mu \in I$. \end{lemma}

{\bf Proof}\\
Suppose that for some $\mu_0$, there exists $\lambda$ such that $f(\lambda)=\mu_0$ and $f'(\lambda)=0$. We should  show that this leads to a contradiction.\\
Let $\lambda_1,\cdots, \lambda_\ell$ the points with this last property. Let $k_j >1$ be the smallest integer such that $f^{(k_j)}(\lambda_j) \neq 0$. Using Rouch\'e's theorem, we see that when $k_j$ is even, necessary $k_j$ complex eigenvalues appear near $\lambda_j$ when  $(\mu -\mu_0) f^{(k_j)}(\lambda_j) <0$ in contradiction with the assumption. Similarly, when $k_j$ is odd, $(k_j-1)$ complex zeros appear when \break  $(\mu -\mu_0) f^{(k_j)}(\lambda_j) \neq 0$.

\begin{lemma}
Let $f(\lambda)$ be a real polynomial of degree $q$ and $g$ a real polynomial of degree $r <q$, such that, for any $\mu \in I=]a,b[$, $f(\lambda)=\mu g(\lambda)$  has $q$ real solutions and suppose that $f$ and $g$ have no common zero,   then $f' g- f g'  \neq 0$, for any $\lambda$ such that $f(\lambda)\in I$. \end{lemma}

{\bf Proof}\\
We have necessarily $g\neq 0$ for these solutions. Hence we can perform the previous argument by applying it to $f/g$.

\begin{proposition}
Except isolated values corresponding to (isolated or embedded) flat bands, the spectrum of the Hou model consists of non overlapping (possibly touching)  bands.
\end{proposition}

Here are  two examples of non trivial closed gaps:
\begin{itemize}
\item
For the triangular model, for $p/q=1/6$, the spectrum is given by :
$$
\{ \lambda\in \mathbb R\,,\, \exists (\theta_1,\theta_2)\in\mathbb R^2,~\lambda^6-18 \lambda^4-12\sqrt{3} \lambda^3+45\lambda^2+36\sqrt{3}\lambda+6-2p^\triangle(6\theta_1,6\theta_2)=0\}$$
i.e. by the condition
$$\lambda^6-18 \lambda^4-12\sqrt{3} \lambda^3+45\lambda^2+36\sqrt{3}\lambda+6  \in[-3,6]\,.
$$
We have
$$
Q_T(\lambda)= \lambda^6-18 \lambda^4-12\sqrt{3} \lambda^3+45\lambda^2+36\sqrt{3}\lambda
$$
which 
satisfies $$Q_T(-\sqrt{3})= Q'_T(-\sqrt{3})= 0\,.$$
 Hence the second gap is closed. Note this is to our knowledge the only closed gap which has been observed for the triangular butterfly (see Figure 2).
\item 
For the graphene model, for $p/q=1/2$, the spectrum is given by
$$
\{ \lambda\in \mathbb R\,,\, \exists (\theta_1,\theta_2)\in\mathbb R^2,~\lambda^4-6\lambda^2+3-2(\cos(2\theta_1)+\cos(2\theta_2)-\cos(2(\theta_1-\theta_2)))\}$$
i.e. 
$$\lambda^4-6\lambda^2\in[-9,0]\,.
$$
 The bands are $[-\sqrt{6},-\sqrt{3}]$, $[-\sqrt{3},0]$, $[0,\sqrt{3}]$ and $[\sqrt{3},\sqrt{6}]$. We have in this case three closed gaps at $-\sqrt{3}, 0,+\sqrt{3}$.
\end{itemize}

\section{Semi-classical analysis for Hou's butterfly near a flat band}\label{s7}

The general study of Hou's butterfly near its flat bands seems difficult, but we can obtain an explicit reduction
for the simplest one, which is the flat band $\{0\}$ in the case when $\omega=0$, $\gamma=4\pi$.
As shown in \cite{KR}, the spectrum of Hou's operator for $\omega=0$, $\gamma=4\pi+h$ is the spectrum of the Weyl $h$-quantization of
\be M(x,\xi,h)=\left(\begin{array}{ccc}
0 & i\,e^{ih/8}(e^{-ix}+e^{-i(x-\xi)}) & -i\,e^{-ih/8}(e^{-ix}+e^{-i\xi}) \\
-i\,e^{-ih/8}(e^{ix}+e^{i(x-\xi)}) & 0 & i\,e^{ih/8}(e^{i(x-\xi)}+e^{-i\xi})\\
i\,e^{ih/8}(e^{ix}+e^{i\xi}) & -i\,e^{-ih/8}(e^{-i(x-\xi)}+e^{i\xi}) & 0
                     \end{array}\right)\ee
 
Let us first recall some rules in semi-classical analysis. The considered  symbols are functions $p(x,\xi,h)$ in the class $S^0(\R^2)$ of smooth functions of $(x,\xi)\in\R^2$ depending on a semi-classical 
parameter $h\in[-h_0,h_0]$, $h_0>0$ (view as ``little'') and satisfying
\be \forall (j,k)\in\N^2\,;\,\exists C_{j,k}\,;\,\forall (x,\xi)\in\R^2,\,
|\partial_x^j\partial_\xi^k p(x,\xi,h)|\leq C_{j,k} \ee
The classical and Weyl quantizations of the symbol $p$ are respectively  (for $h\neq 0$, $|h|\leq h_0$) the pseudodifferential operators acting on $L^2(\R)$ by 
\ba p(x,hD_x,h)u(x)&=&\frac{1}{2\pi h}\int\!\!\!\int e^{i(x-y)\xi/h} p(x,\xi,h)\, u(y)\,dy\,d\xi \,,\\
\text{Op}^W_h(p)u(x)&=&\frac{1}{2\pi h}\int\!\!\!\int e^{i(x-y)\xi/h} p(\frac{x+y}{2},\xi,h)\, u(y)\,dy\,d\xi \,.\ea
 Conversely, if $P$ is a pseudodifferential operator, we denote $\sigma(P)$ and $\sigma^W(P)$ its classical and Weyl symbols. If
these symbols admit asymptotic expansions 
\ban \sigma(P)(x,\xi,h)=\sigma_0(P)(x,\xi)+h\,\sigma_{-1}(P)(x,\xi)+\mathcal O(h^2)\,,\\
\sigma^W(P)(x,\xi,h)=\sigma_0^W(P)(x,\xi)+h\,\sigma_{-1}^W(P)(x,\xi)+\mathcal O(h^2)\,\ean
they are related by
\ba \sigma_0^W(P)(x,\xi)&=&\sigma_0(P)(x,\xi)\,, \label{prin}\\
\sigma_{-1}^W(P)(x,\xi)&=&\sigma_{-1}(P)(x,\xi)-\frac{1}{2i} \partial_x\partial_\xi \sigma_0(P)(x,\xi)\,. \label{subprin}\ea
$\sigma_0^W(P)$ and $\sigma_{-1}^W(P)$ are called the principal and subprincipal symbols of $P$. If $P$ and $Q$ are pseudodifferential
operators admitting such expansions, the classical composition\footnote{By this, we mean that we use the pseudo-differential calculus involving the classical quantization.} is given by
\be \sigma_0(P\,Q)=\sigma_0(P)\,\sigma_0(Q)\,,~~\sigma_{-1}(P\,Q)=\sigma_{-1}(P)\,\sigma_0(Q)+\sigma_0(P)\,\sigma_{-1}(Q)+\frac{1}{i} \partial_\xi P\,\partial_x Q \label{comp} \ee
Another important fact, which partially justifies the use of Weyl quantization in the study of selfadjoint operators, is
\be \sigma^W(P^*)=\left(\sigma^W(P)\right)^* \label{selfadjoint}\,.\ee

In our case, the principal symbol $M_0$ is given by 
\be M_0(x,\xi)=\left(\begin{array}{ccc}
0 & i(e^{-ix}+e^{-i(x-\xi)}) & -i(e^{-ix}+e^{-i\xi}) \\
-i(e^{ix}+e^{i(x-\xi)}) & 0 & i(e^{i(x-\xi)}+e^{-i\xi})\\
i(e^{ix}+e^{i\xi}) & -i(e^{-i(x-\xi)}+e^{i\xi}) & 0
                     \end{array}\right)\ee

We first prove :

\begin{proposition}
 There exists a familly $U_0(x,\xi)$ of unitary $3\times 3$ matrices, depending smoothly on $(x,\xi)$, $2\pi$-periodic 
in each variable, and a familly $A(x,\xi)$ of selfadjoint $2\times 2$ matrices such that
\be U_0^*(x,\xi)\,M_0(x,\xi)\,U_0(x,\xi)= \left(\begin{array}{lc}
                                            ~\,0 & \begin{array}{cc}
0 & 0 
                                           \end{array}\\
\begin{array}{l}
 0 \\0
\end{array} & A(x,\xi)\end{array}\right)\,.\ee
Moreover,  for any $(x,\xi)\in\R^2$, the spectrum of $A(x,\xi)$ is contained in $[-2\sqrt{3},-\sqrt{3}]\cup[\sqrt{3},2\sqrt{3}]$.
\end{proposition}

{\bf Proof :} We easily compute the characteristic polynomial 
\be \det(M_0(x,\xi)-\lambda\,I_3)=-\lambda^3+(6+2p^\triangle(x,\xi))\lambda \,.\ee
The range of $p^\triangle$ is $[-3/2,3]$, so the kernel of $M_0(x,\xi)$ has dimension 1, and
the spectrum of the restriction of $M_0$ to $(\ker(M_0(x,\xi))^\perp$ is contained in $[-2\sqrt{3},-\sqrt{3}]\cup[\sqrt{3},2\sqrt{3}]$.
A unitary basis vector of $\ker(M_0(x,\xi))$ is $e_0(x,\xi)=\alpha(x,\xi)\,\tilde e_0(x,\xi)$ with
\be \tilde e(x,\xi)=\left(\begin{array}{c}
                           1+e^{-ix} \\
1+e^{i(x-\xi)} \\
1+e^{i\xi}
                          \end{array}\right)\ee
\be \alpha(x,\xi)=\frac{1}{\sqrt{6+2p^\triangle (x,\xi)}} \ee
So we choose $e_0(x,\xi)$ as the first column of $U_0(x,\xi)$. We then observe
\be \rm{Re}\left\langle\tilde e_0(x,\xi),\left(\begin{array}{l}
                          1 \\ 1 \\ 1
                         \end{array}\right)\right\rangle=3+p^\triangle(x,\xi)\geq\frac{3}{2}\,, \ee
and thus consider a unitary $3\times 3$ matrix $B$ whose first line is $\frac{1}{\sqrt{3}} (1,1,1)$. Then
\be B\,e_0(x\xi)=\left(\begin{array}{c}
                        a(x,\xi) \\ b(x,\xi)\\ c(x,\xi)
                       \end{array}\right) \ee
where $\rm{Re}(a(x,\xi))>0$. We define the unitary vector $f(x,\xi)$ by
\be f(x,\xi)=B^*\frac{1}{\sqrt{|a(x,\xi)|^2+|b(x,\xi)|^2}}\left(\begin{array}{c}
                                                             -\bar b(x,\xi) \\ \bar a(x,\xi) \\ 0
                                                            \end{array}\right)\ee
$f(x,\xi)$ is orthogonal to $e_0(x,\xi)$ and we put
\be g(x,\xi)=\overline{e_0(x,\xi)\wedge f(x,\xi)} \,.\ee
We finally take $U_0(x,\xi)=(e_0(x,\xi),f(x,\xi),g(x,\xi))$.\\ \qed\\
 \begin{remark}
We have preferred to give a complete elementary proof for the triviality of the fiber bundle whose fiber at $(x,\xi)$ is the eigenspace of $M(x,\xi)$ associated with  the two non vanishing eigenvalues. As observed by G. Panati, this can be obtained by general results (see in particular Proposition 4  in \cite{Pan}).
\end{remark}

Using Proposition 3.3.1 in \cite{HS2} and its corollary, we get:
\begin{proposition}
 There exist a unitary $3\times 3$ pseudodifferential operator $U$ with principal symbol $U_0(x,\xi)$, 
a selfadjoint scalar operator $\mu$ with principal symbol $0$, and a selfadjoint $2\times 2$ operator $\tilde A$
with principal symbol $A(x,\xi)$ such that
\be U^*\,\rm{Op}_h^W(M(x,\xi,h)) \,U=\left(\begin{array}{lc}
                                            \;\mu & \begin{array}{cc}
0 & 0 
                                           \end{array}\\
\begin{array}{l}
 0 \\0
\end{array} & \tilde A \end{array}\right)\ee
Moreover, the part of the spectrum of $\rm{Op}_h^W(M(x,\xi,h))$ in any compact subset of $]-\sqrt{3},\sqrt{3}[$ is that one of $\mu$ 
for $|h|$ small enough.
\end{proposition}

The main result of this section is the computation of the
subprincipal symbol of $\mu$.

\begin{proposition} : 
\be \sigma^W(\mu)(x,\xi,h)=-h\frac{3-p^\triangle(x,\xi)}{4(3+p^\triangle(x,\xi))}+\mathcal O(h^2)\,. \ee
\end{proposition}

{\bf Proof :} The computation is in the spirit of \S 6.2 in \cite{HS2}. In this text\footnote{ Note that one term has disappeared at the printing in formula (6.2.9) in \cite{HS2} which is fortunately re-established in formula (6.2.19).}   the matrix
$M(x,\xi)$ satisfies in addition  $\partial_x \partial_\xi M(x,\xi)=0$ and does not depend on $h$. On the other hand, we are here helped by
the relation $M_0(x,\xi)e_0(x,\xi)=0\,$.\\

Since $\sigma_0(\mu)=0$, (\ref{subprin}) gives 
\be \sigma_{-1}^W(\mu)=\sigma_{-1}(\mu)=\sigma_{-1}(U^*\,\rm{Op}_h^W(M(x,\xi,h))\, U)_{11}\,.\ee 
We use the classical calculus to compute this term. Let $U(x,\xi,h)$, 
$V(x,\xi,h)$ and 
$$ N(x,\xi,h)=N_0(x,\xi)+h\,N_1(x,\xi)+\mathcal O(h^2)$$ be the classical symbols of $U$, $U^*$ and $\rm{Op}_h^W(M(x,\xi,h))$.\\
 Using (\ref{prin}), (\ref{subprin}) and (\ref{selfadjoint}) we observe~:
\begin{enumerate}
\item The first column of $U(x,\xi,h)$ is on the form $e_0(x,\xi)+h\,e_1(x,\xi)+\mathcal O(h^2)$\,.
\item The first line of $V(x,\xi)$ is on the form $\bar e_0^T(x,\xi)+h\,f_1^T(x,\xi)+\mathcal O(h^2)$. 
\item $N_0(x,\xi)=M_0(x,\xi)$\,.
\item $N_1(x,\xi)=M_1(x,\xi)+\frac{1}{2i}\partial_x \partial_\xi M_0(x,\xi)\,$.
\end{enumerate}
Then the composition rules (\ref{comp}) together with $$M_0(x,\xi)\,e_0(x,\xi)=0$$ and $$\bar e_0^T(x,\xi)\,M_0(x,\xi)=0$$ give~:
\ban 
\sigma^W_{-1}(\mu)(x,\xi)&=&f_1^T(x,\xi) \,N_0(x,\xi) \,e_0(x,\xi) + \bar e_0^T(x,\xi)\, N_1(x,\xi) \,e_0(x,\xi)+ 
\bar e_0^T(x,\xi)\, N_0(x,\xi)\, e_1(x,\xi)\\
&&+\frac{1}{i}\partial_{\xi}(\bar e_0^T(x,\xi)\, N_0(x,\xi))\,\partial_x e_0(x,\xi)
+(\frac{1}{i}\partial_{\xi}\bar e_0^T(x,\xi) \partial_x N_0(x,\xi))\, e_0(x,\xi) \\
&=& \bar e_0^T(x,\xi) \,M_1(x,\xi)\, e_0(x,\xi)+ 
\frac{1}{2i}\bar e_0^T(x,\xi)\,\partial_x\partial_{\xi}M_0(x,\xi)\, e_0(x,\xi)\\
&&+\frac{1}{i}(\partial_{\xi}\bar e_0^T(x,\xi) \,\partial_x M_0(x,\xi)\, e_0(x,\xi)) \,. \ean
Then differentiating the identity  $M_0(x,\xi)e_0(x,\xi)=0$ successively gives: 
\ban \partial_x M_0(x,\xi)\,e_0(x,\xi)&=&-M_0(x,\xi)\,\partial_x e_0(x,\xi)\,,\\
\partial_{\xi} M_0(x,\xi)\,e_0(x,\xi)&=&-M_0(x,\xi)\,\partial_{\xi} e_0(x,\xi)\,,\\
\bar \partial_x \partial_{\xi} M_0(x,\xi) \,e_0(x,\xi)&=&
-\partial_x M_0(x,\xi)\,\partial_{\xi} e_0(x,\xi)
-\partial_{\xi} M_0(x,\xi)\,\partial_{x} e_0(x,\xi)\\
&&-M_0(x,\xi)\,\partial_x \partial_\xi e_0(x,\xi) \,,\\
\bar e_0^T(x,\xi)\, \partial_x \partial_{\xi} M_0(x,\xi) \,e_0(x,\xi)&=&
\partial_x \bar e_0^T(x,\xi)\, M_0(x,\xi)\,\partial_{\xi} e_0(x,\xi)\\
&&+\partial_{\xi} \bar e_0^T(x,\xi)\, M_0(x,\xi)\,\partial_{x} e_0(x,\xi) \,.\ean
Hence
\be \sigma^W_{-1}(\mu)(x,\xi)=\langle M_1(x,\xi)\, e_0(x,\xi),e_0(x,\xi)\rangle 
+\rm{Im}\langle M_0(x,\xi)\,\partial_{\xi}e_0(x,\xi),\partial_x e_0(x,\xi)\rangle \,.\ee
A straightforward computation gives
\be \rm{Im}\langle M_0(x,\xi)\,\partial_{\xi}e_0(x,\xi),\partial_xe_0(x,\xi)\rangle = 
\frac{p^\triangle(x,\xi)}{3+p^\triangle(x,\xi)} \,.\ee
On the other side, we denote by $\lambda(x,\xi,h)$ the second eigenvalue of $M(x,\xi,h)$. The computation of the characteristic polynomial $\det(M(x,\xi)-\lambda \,I_3)$ gives
\be -\lambda^3(x,\xi,h)+(6+2p^\triangle(x,\xi))\lambda(x,\xi,h) +4\sin \frac{3h}{8}(1+p^\triangle(x,\xi))=0 \,.\ee
So
\ba \langle M_1(x,\xi)\, e_0(x,\xi),e_0(x,\xi)\rangle &=&\langle \partial_h M(x,\xi,0) \,e_0(x,\xi),e_0(x,\xi)\rangle \\
&=&\partial_h \lambda(x,\xi,0)=-\frac{3(1+p^\triangle(x,\xi))}{4(3+p^\triangle(x,\xi))}\,.
\ea
Hence 
\ban \sigma^W_{-1}(\mu)(x,\xi)
&=&-\frac{3(1+p^\triangle(x,\xi))}{4(3+p^\triangle(x,\xi))} + \frac{p^\triangle(x,\xi)}{3+p^\triangle(x,\xi)} \\
&=& -\frac{3-p^\triangle(x,\xi)}{4(3+p^\triangle(x,\xi))} \,.\ean
Then $\sigma^W_0(\mu)=0$ achieves the proof.\\ \qed \\
 \section{Conclusion}
In this paper we have shown that for the model proposed by Hou relative to the kagome lattice and whose justification for the analysis of the Schr\"odinger magnetic operator was given in \cite{KR}, a Chambers analysis is available permitting to recover most of the characteristics observed in the case of the square lattice for the Hofstadter butterfly, the triangular butterfly  or the hexagonal (graphene) butterfly. This makes all the semi-classical techniques developed in \cite{HS1,HS2,Ke} available but this leads also to new questions to analyze: the existence of flat bands. In the previous section we have shown how, when the flux is close to $4\pi$ ($\gamma =4\pi +h$) the semi-classical calculus permits via the computation of a subprincipal symbol to reduce the spectral analysis of the Hou operator in the interval $[-\sqrt{3} + \epsilon_0, \sqrt{3} -\epsilon_0]$ ($\epsilon_0>0$) to the analysis of a $h$- pseudodifferential operator with explicit principal symbol. In particular, our analysis implies that the convex hull of the part of the spectrum contained in this interval is
 $[-\frac 34 h + \mathcal O (h^2)\,,\, \mathcal O (h^2)]$ for $h>0$ and $[ \mathcal O (h^2)\,,\,  -\frac 34 h +\mathcal O (h^2)]$ 
  for $h<0$. This suggests the beginning of a renormalization involving after one step the perturbation of a function of the triangular Harper model. More precisely, this function is the function
  $$
  \lambda \mapsto  -h \frac{3-\lambda}{4 (3+\lambda)}\,.
  $$

{\bf Acknowledgements.}\\
B. Helffer and P. Kerdelhu\'e are partially supported by the ANR programme NOSEVOL.  B. Helffer thanks Gianluca Panati for useful discussions and Gian Michele Graf for the transmission of \cite{AEG}.\\

\begin{figure}[ht!]
\begin{center}
\includegraphics[width=9cm]{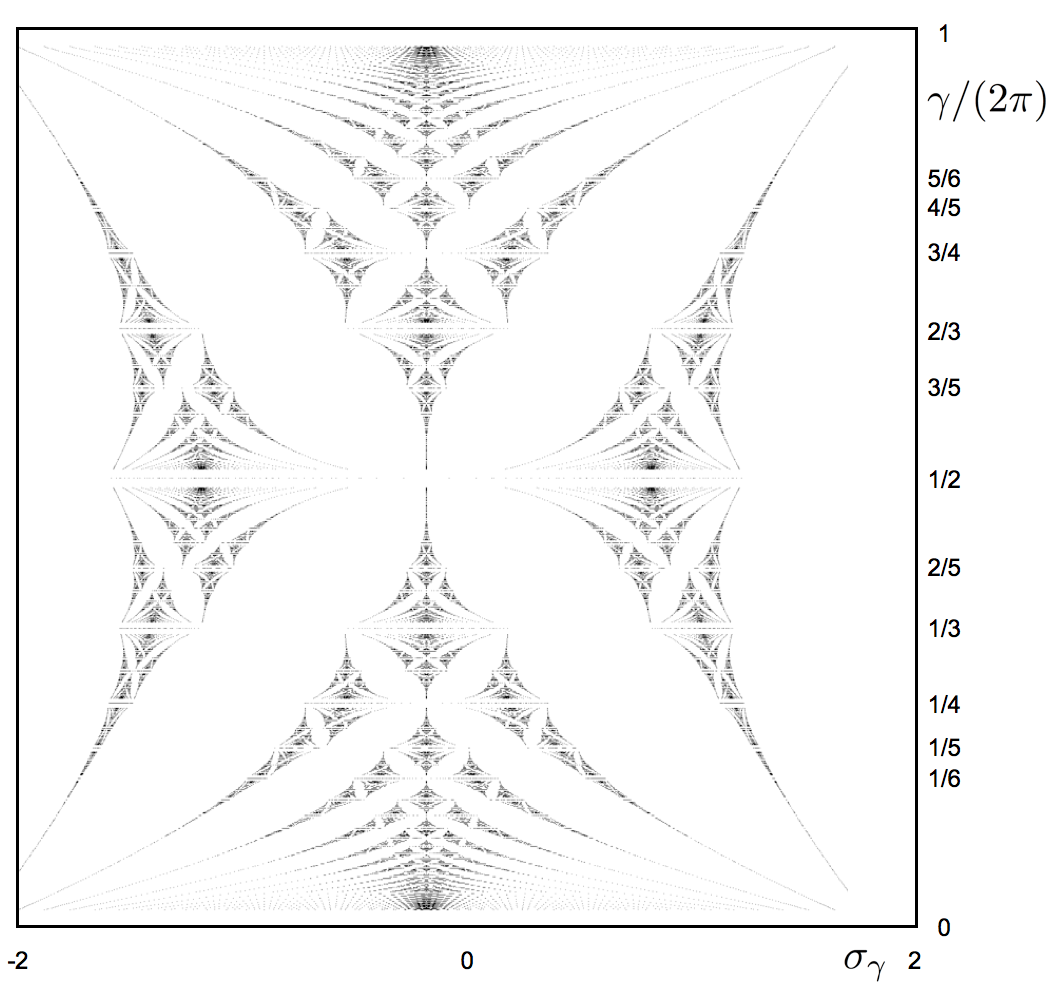}
\caption{Square lattice.}
 \end{center}
 \end{figure}
 
 \begin{figure}[ht!]
\begin{center}
\includegraphics[width=10cm]{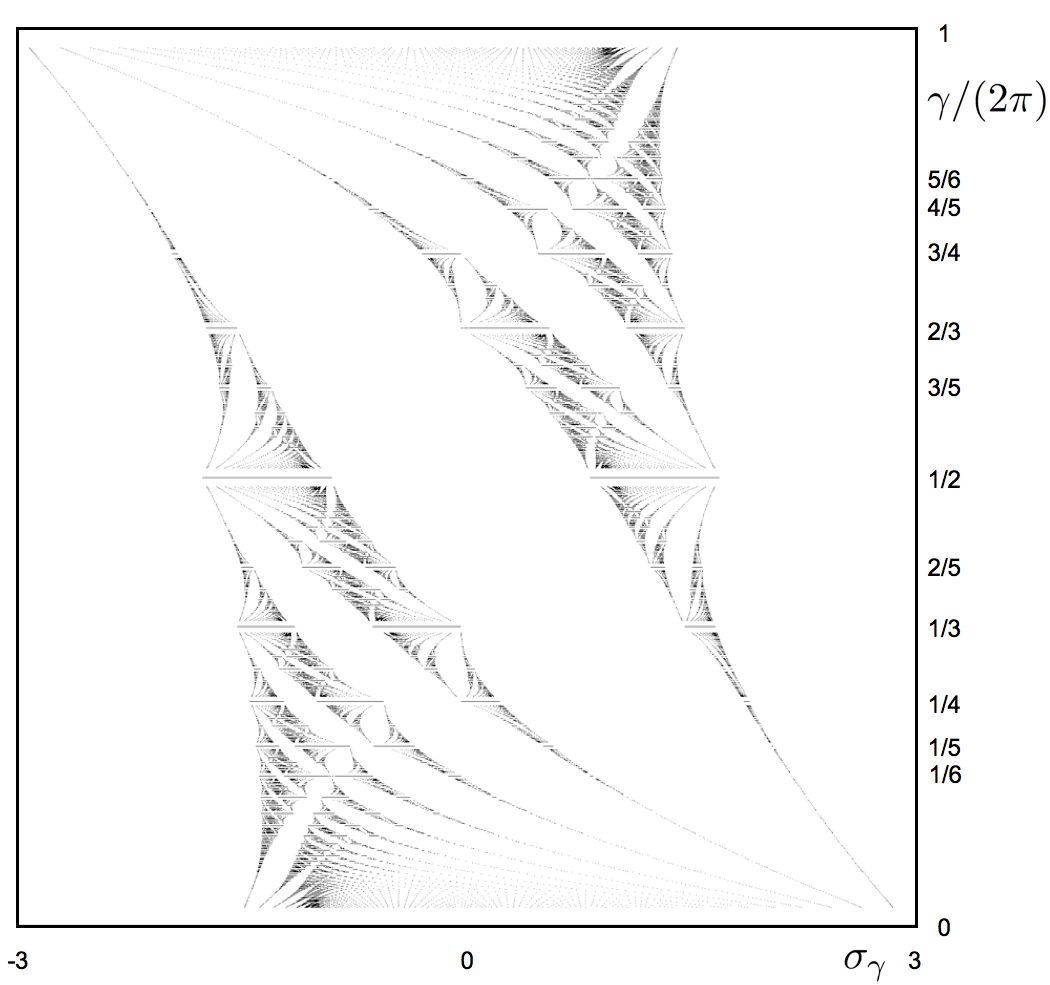}
\caption{triangular lattice.}
 \end{center}
 \end{figure}

 \begin{figure}[ht!]
\begin{center}
\includegraphics[width=10cm]{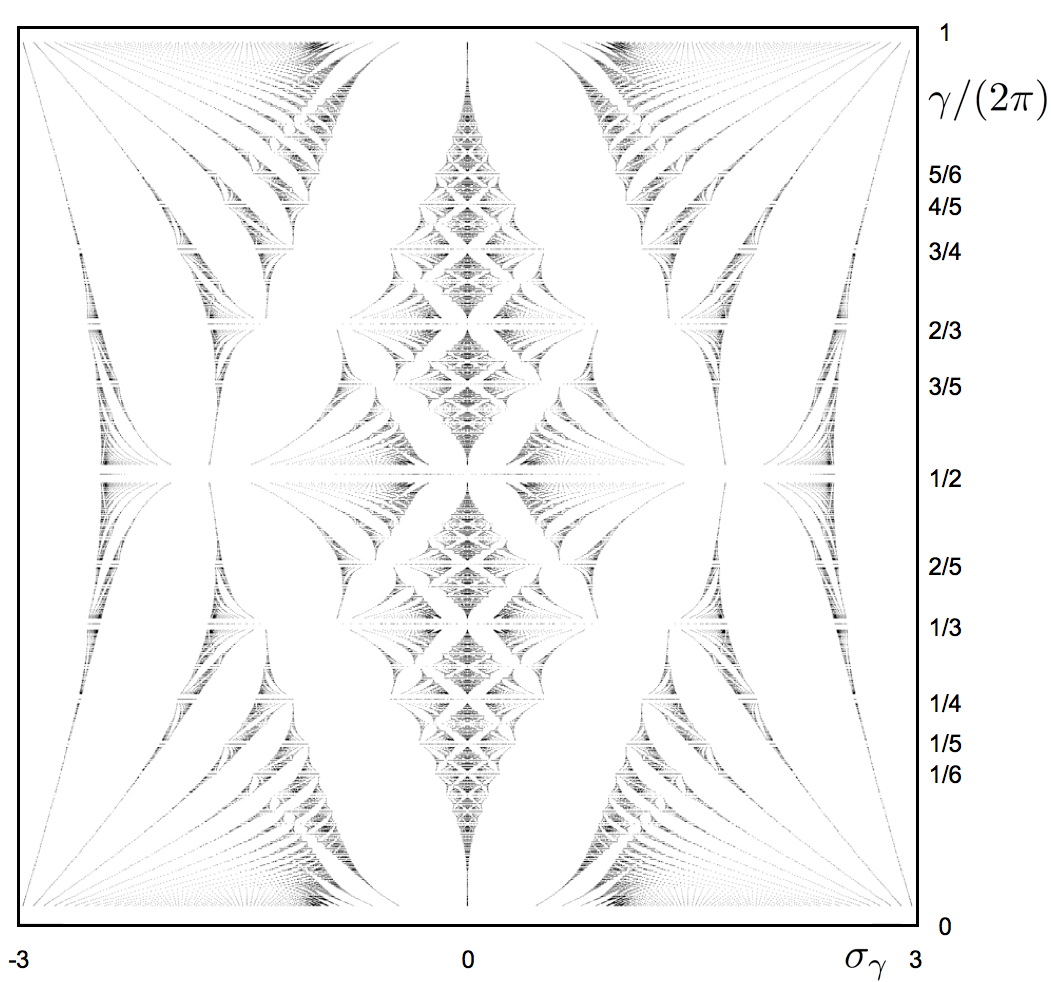}
\caption{Hexagonal  lattice.}
 \end{center}
 \end{figure}
 
  \begin{figure}[ht!]
\begin{center}
\includegraphics[width=14cm]{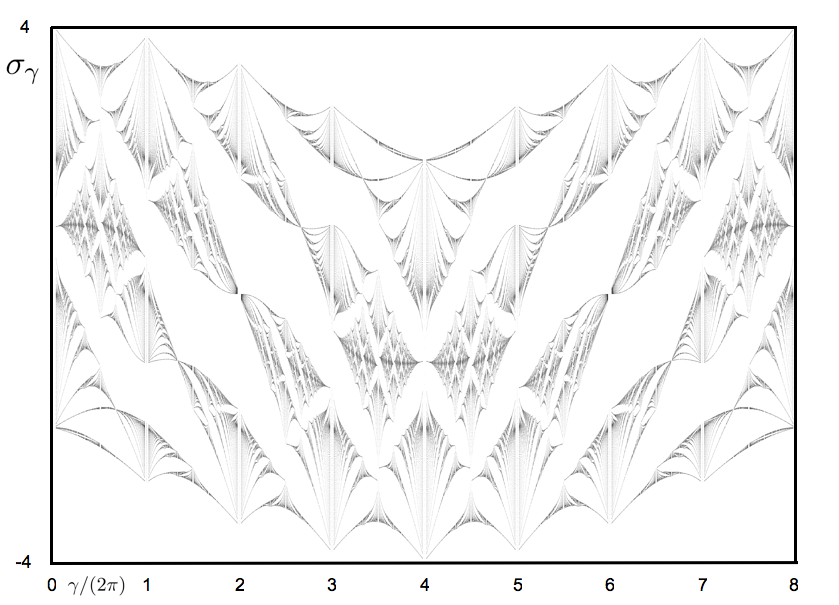}
\caption{Kagome lattice, $\omega=0$.}
 \end{center}
 \end{figure}
 
   \begin{figure}[ht!]
\begin{center}
\includegraphics[width=14cm]{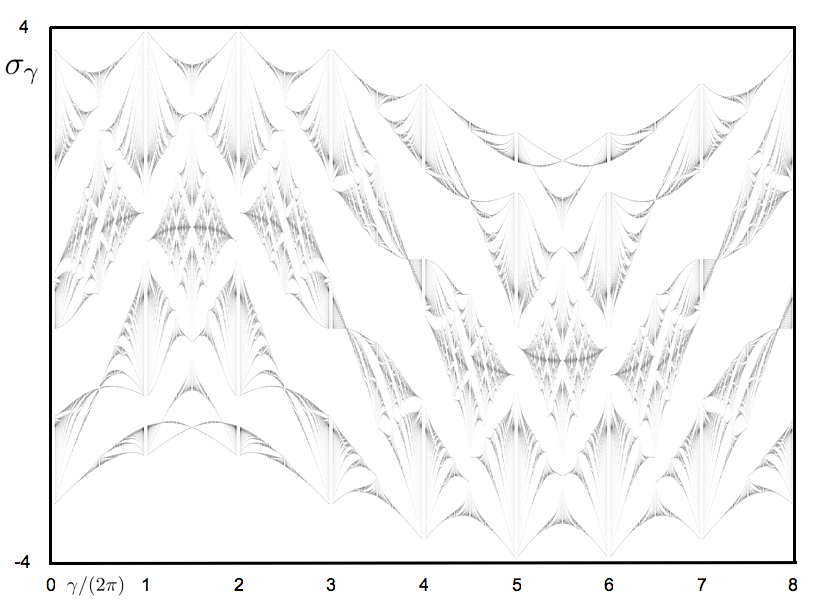}
\caption{Kagome lattice, $\omega=\frac \pi 8$.}
 \end{center}
 \end{figure}


\begin{thebibliography}{99}
\bibitem{AEG} A. Agazzi, J.-P. Eckmann, and G.M. Graf.
\newblock The colored Hofstadter butterfly for the Honeycomb lattice.
\newblock arXiv:1403.1270v1 [math-ph] 5 Mar 2014.

\bibitem{AJ} A. Avila and S. Jitomirskaya. 
\newblock The Ten Martini Problem. 
\newblock Ann. of Math., 170 (2009) 303--342.

\bibitem{AKY} J. E. Avron, O. Kenneth and G. Yeshoshua.
\newblock A numerical study of the window condition for
Chern numbers of Hofstadter butterflies.
\newblock arXiv:1308.3334v1 [math-ph], 15 Aug 2013.

\bibitem{Az} Ya. Azbel.
\newblock Energy spectrum of a conduction electron in a magnetic field.
\newblock  Sov. Phys.  JETP 19 (3) (1964),  264.

\bibitem  {bel} J. Bellissard. 
\newblock Le papillon de Hofstadter.
\newblock  Ast\'erisque {\bf 206} (1992) 7--39.

\bibitem{BKS} J. Bellissard, C. Kreft, and  R. Seiler.
\newblock Analysis of the spectrum of a particle on a triangular lattice with two magnetic fluxes by algebraic  and numerical methods.
\newblock J.~of Physics A24 (1991),  2329-2353.

\bibitem{belsim} J. Bellissard and B. Simon.
\newblock   Cantor spectrum for the almost Mathieu equation.
\newblock J.~Funct. Anal. 48,  408-419  (1982).

\bibitem{Ch} W. Chambers.
\newblock Linear network model for magnetic breakdown in two dimensions. 
\newblock Phys. Rev A140 (1965), 135--143.

\bibitem {CW} F.H. Claro and G.H.  Wannier.
\newblock Magnetic subband structure of electron in hexagonal
lattices. 
\newblock Phys. Rev. B, Volume 19, No 12 (1979),  6068--6074.

\bibitem{DM} P. Delplace and G. Montambaux.
\newblock WKB analysis of edge states in graphene in a strong magnetic field.
\newblock arXiv:1007.2910v1 [cond-mat.mes-hall] 17 Jul 2010.

\bibitem{Ha} P.G. Harper.
\newblock Single band motion of conduction electrons in a uniform magnetic field.
\newblock  Proc. Phys. Soc. London A 88 (1955), 874.

\bibitem {HS1}
B. Helffer and J. Sj\"ostrand.
\newblock Analyse semi-classique pour l'\'equation de Harper (avec application \`a l'\'equation de Schr\"odinger avec 
champ magn\'etique).
M\'em. Soc. Math. France (N.S.)  {\bf 34}  (1988) 1--113.

\bibitem  {HS2}
B. Helffer and J. Sj\"ostrand. \newblock Analyse semi-classique pour l'\'equation de Harper. II.
Comportement semi-classique pr\`es d'un rationnel.
\newblock  M\'em. Soc. Math. France (N.S.) 
{\bf 40}  (1990) 1--139.


\bibitem  {Hof}
D. Hofstadter.  Energy levels and wave functions of Bloch electrons in rational and irrational magnetic fields.
Phys. Rev. B {\bf 14} (1976) 2239--2249.

\bibitem  {Hou} J-M. Hou.
\newblock Light-induced Hofstadter's butterfly spectrum of ultracold atoms on the two-dimensional kagome lattice, 
\newblock CHN. Phys. Lett. {\bf 26}, 12 (2009), 123701.


\bibitem{Ke} P. Kerdelhu\'e.
\newblock Spectre de l'op\'erateur de Schr\"odinger magn\'etique avec
sym\'etrie
 d'ordre 6.
\newblock M\'emoire de la SMF, tome 51 (1992),  1--139.


\bibitem{KR} P. Kerdelhu\'e and J. Royo-Letelier.
\newblock On the low lying spectrum of the magnetic Schr\"odinger operator with kagome periodicity.
\newblock arXiv:1404.0642v1 [math.AP] 2 April 2014. (submitted).

\bibitem {Kr} C. Kreft.
\newblock Spectral analysis of Hofstadter-like models.
\newblock Thesis Technischen Universit\"at Berlin, Berlin (1995).

\bibitem{Pan} G. Panati.
\newblock Triviality of Bloch and Bloch-Dirac bundles.
\newblock Annales Henri Poincar\'e 8 (5), 995-1011.

\bibitem{ReSi} M. Reed and B. Simon.
\newblock  Methods of modern mathematical physics. {IV}.
\newblock  Analysis of operators, Academic Press, New York, 1978.


\bibitem{VM} P. Van Mouche.
\newblock The coexistence problem for the discrete Mathieu operator.
Comm. Math. Phys. 122 (1989), no. 1, 23--33.


\bibitem {WiAu} M. Wilkinson and E. Austin.
\newblock Semi-classical analysis of phase space lattices with three
fold
 symmetry.
\newblock J. Phys. A: Math. Gen. 23,  2529--2553.   

\end{thebibliography}
 \end{document}